\documentclass[12pt]{article}

\usepackage{amsmath,amssymb,amsfonts,theorem,makeidx,latexsym,epsfig,subfigure}

\usepackage{color}

\newtheorem{defn}{Definition}[section]

\newtheorem{lemma}[defn]{Lemma}

{\theorembodyfont{\rmfamily}

\newtheorem{ex}[defn]{Example}}

\newtheorem{thm}[defn]{Theorem}

\newtheorem{prop}[defn]{Proposition}

\newtheorem{cor}[defn]{Corollary}

\newtheorem{rem}[defn]{Remark}

\numberwithin{equation}{section}

\newcommand{\h}{{\cal H}}

\newcommand{\ltr}{ L^2(\mathbb R) }

\newcommand{\ltn}{{\ell}^2(\mathbb N)}

\newcommand{\si}{S^{-1}}

\newcommand{\mn}{\mathbb N}

\newcommand{\mr}{\mathbb R}

\newcommand{\mz}{\mathbb Z}

\newcommand{\mc}{\mathbb C}

\newcommand{\mts}{ \{E_{mb}T_{na}g \}_{m,n \in \mz}}

\def\bp{{\noindent\bf Proof. \ }}

\def\ep{\hfill$\square$\par\bigskip}

\def\bqs{\begin{equation}}

\def\eqs{\tag*{$\square$}\end{equation}\par\bigskip}

\def\hp{\hat{\psi}}

\def\la{\langle}

\def\ra{\rangle}

\def\ga{\gamma}

\def\ctk{\{c_k\}_{k=1}^\infty}

\def\etk{\{e_k\}_{k=1}^\infty}

\def\suk{\sum_{k=1}^\infty}

\def\sujz{\sum_{j\in \mz}}

\def\hp{\widehat{\psi}}

\def\span{\overline{\text{span}}}

\def\Span{\text{span}}

\def\vn{\vspace{.1in}\noindent}

\def\bop{\begin{op}\rm}

\def\eop{\end{op}}

\def\bee{\begin{eqnarray}}

\def\ene{\end{eqnarray}}

\def\bes{\begin{eqnarray*}}

\def\ens{\end{eqnarray*}}

\def\bei{\begin{itemize}}

\def\eni{\end{itemize}}

\def\bt{\begin{thm}}

\def\et{\end{thm}}

\def\bc{\begin{cor}}

\def\ec{\end{cor}}

\def\bpr{\begin{prop}}

\def\epr{\end{prop}}

\def\bl{\begin{lemma}}

\def\el{\end{lemma}}

\def\bd{\begin{defn}}

\def\ed{\end{defn}}

\def\bex{\begin{ex}}

\def\enx{\end{ex}}

\def\bfi{\begin{fig}}

\def\efi{\end{fig}}

\def\inr{\int_{-\infty}^\infty}

\newcommand{\nft}{ || f||^2}

\def\hp{\widehat{\psi}}

\def\wfi{\{D_{a^j}T_{kb}\psi\}_{j,k\in \mz}}

\def\wfN{\{D_{a^{Nj}}T_{kb}\psi\}_{j,k\in \mz}}

\def\wfb{\{D_{a^j}T_{kbN}\psi\}_{j,k\in \mz}}

\def\wtail{\sum_{k\neq 0} \sujz |\hp(a^j\ga) \hp(a^j\ga +k/b)|   }

\title{Weaving Information Packets\thanks{The research was supported by
the Basic Science Research Program through the National Research
Foundation of Korea (NRF) funded by the Ministry of Education
(2016R1D1A1B02009954).}}

\date{\today}

\author{Ole Christensen, Hong Oh Kim, Rae Young Kim}

\begin{document}

\maketitle

\begin{abstract} The concept of weaving
     of frames for Hilbert spaces
     was introduced by Bemrose et al. in 2016.
     Two frames $\{f_k\}_{k\in I},
    \{g_k\}_{k\in I}$ are  woven if the ``mixed system"
    $\{f_k\}_{k\in \sigma} \cup \{g_k\}_{k\in I\setminus \sigma}$
    is a frame for each index set $\sigma \subset I;$ that is, processing a signal using two woven frames yields a certain stability against loss of information. The concept
    easily extends to $N$ frames, for any integer $N>2.$  Unfortunately it is nontrivial to construct  useful woven frames, and the literature is sparse concerning explicit constructions. In this paper we introduce so-called information packets, which contain as well frames as fusion frames as special case. The concept of woven frames immediately generalizes to information packets, and we demonstrate how to construct practically relevant   woven information packets based on particular wavelet systems in $\ltr.$ Interestingly, we show that certain wavelet systems can be split into $N$ woven information packets, for any integer $N\ge 2.$ We finally consider corresponding questions for  Gabor system in $\ltr,$ and prove that for any fixed $N\in \mn$ we can find a Gabor frame that can be split into $N$ woven information packets; however,
    in contrast to the wavelet case, the density conditions for Gabor system excludes the
    possibility of finding a single Gabor frame
    that works simultaneously for all $N\in \mn.$
    \end{abstract}

\begin{minipage}{120mm}


\end{minipage}
\

\section{Introduction}

 Our work is inspired by the concept of  weaving frames, introduced in \cite{CasBem}. Two frames $\{f_k\}_{k\in I},
 \{g_k\}_{k\in I}$ for a Hilbert space $\h$  are said to be {\it woven} if the ``mixed system"
 $\{f_k\}_{k\in \sigma} \cup \{g_k\}_{k\in I\setminus \sigma}$
 is a frame for each index set $\sigma \subset I.$ Having a set of two woven frames at hand yields a particular type of robustness towards loss of information, because the frame property is kept even if some of the elements $f_k$ are erased (as long as the corresponding elements $g_k$ are kept).

 As is clear already from \cite{CasBem}, it is notoriously difficult to construct interesting cases of woven frames; see also the recent papers \cite{CKKM25,CabMol}. The purpose
 of the current paper is to introduce the idea of weaving in  the more general context of what we will call {\it information packets} and to provide explicit examples using wavelet systems in $\ltr.$

 The technical definition of information packets and weaving hereof will be given in Section \ref{251701a}, but let us end this introduction with a short motivation for this concept.
 Let $\h$ denote a separable Hilbert space.
 Several variants of frame theory deals with local analysis of
 signals $f\in \h.$ While classical frame theory provides a global
 decomposition of signals $f\in \h,$ {\it fusion frames} consider
 first a decomposition of $\h$ into a number of subspaces $W_j, \, j\in J,$ with the purpose of providing a {\it local decomposition}
 within each space and then piece them together to obtain a decomposition of $\h.$  Information packets  contain fusion frames as a special case, and are more suitable for our purpose of considering weaving problems.

 The paper is organized as follows. In the rest of this section we provide the necessary background information about fusion frames and
 wavelet frames. In Section \ref{251701a} we state the formal definition of information packets and relate it to other types of decomposition of Hilbert spaces. Finally we introduce weaving of information packets in Section \ref{251901s}   and provide an explicit construction of woven information packets based on wavelet systems in $\ltr$.
 More precisely,
 we show that certain wavelet systems can be split into $N$ woven information packets, for any integer $N\ge 2.$
 We conclude the paper with a short discussion of weaving in the context of Gabor systems in $\ltr,$ with special focus on the differences between the wavelet case and the Gabor case. In particular, while for each $N\ge 2$ it is indeed possible to construct a Gabor frame that can be split into $N$ woven information packets, the density condition for Gabor frames excludes that a single Gabor frame can be split into $N$ woven information packets for all $N\ge 2.$

\subsection{Fusion frames}
Fusion frames were introduced in the paper \cite{CasKuLi}, as follows:

\bd \label{251801c} Consider a sequence $\{\omega_j\}_{j\in J}$ of strictly positive scalars and a sequence  $\{W_j\}_{j\in J}$ of closed subspaces of a Hilbert space $\h.$ Denoting the orthogonal projection of $\h$ onto $W_j$ by $P_j,$ \,  $\{W_j, \omega_j\}_{j\in J}$ is called a fusion frame if there exist constants $A,B>0$ such that
\bes A\, \nft \le \sum_{j\in J} \omega_j^2 || P_jf||^2 \le B\, \nft, \, \forall f\in \h.\ens  \ed

It is fair to say that the theory for fusion frames in most aspects runs precisely like classical frame theory, with results and proofs that are closely related. For instance,
one can prove \cite{CasKuLi} that if  $\{W_j, \omega_j\}_{j\in J}$ is a fusion frame, the {\it fusion frame operator}
\bes S: \h \to \h, Sf:= \sum_{j\in J} \omega_j^2  P_jf\ens is
bounded and bijective, furthermore, the infinite series defining $S$ is unconditionally convergent for all $f\in \h.$ This immediately  leads to the
{\it fusion frame decomposition}
\bee \label{241901e} f=  \sum_{j\in J} \omega_j^2 S^{-1} P_jf, \, f\in \h,\ene
again with unconditional convergence.

\subsection{Wavelet frames $\wfi$} \label{242507c}

A (standard) wavelet system in $\ltr$ is a collection of functions having the form \\
$\{a^{j/2}\psi(a^jx-kb)\}_{j,k\in \mz}$ for some $a>1, b>0,$ and
$\psi \in \ltr.$ By introducing the scaling operator
$D_a f(x):= a^{1/2}f(ax)$ and the translation operator $T_bf(x):=f(x-b),$ a wavelet system can be written as $\wfi.$

There is a vast literature on wavelet frames, see, e.g., \cite{Da2,CB}. For our analysis, we need the following
result, first stated in \cite{CC3}  (see also  Theorem 19.1.1 in \cite{CB}).  It deals
with more general wavelet-type systems of the form
 $\{D_{\lambda_j}T_{kb} \psi \}_{j,k\in\mz},$  where
$\{\lambda_j\}_{j\in\mz}$ is any  sequence of positive real scalars.

\bl \label{r-1}
Let $\{\lambda_j\}_{j\in\mz}$ be a sequence of positive real scalars, $b>0$ and
$\psi \in \ltr$.
Suppose that
\bee \label{r-1b}
B:=\frac1{b} \sup_{\ga \in \mr} \left( \sujz
| \hp(\frac{\ga}{\lambda_j})|^2  + \sum_{k\neq 0} \sujz  | \hp(\frac{\ga}{\lambda_j})
\hp(\frac{\ga}{\lambda_j} +\frac{k}{b})| \right) < \infty
\ene and

\bee \label{r-1a}
A:=\frac1{b} \inf_{\ga \in \mr} \left( \sujz | \hp(\frac{\ga}{\lambda_j})|^2
- \sum_{k\neq 0} \sujz |\hp(\frac{\ga}{\lambda_j}) \hp(\frac{\ga}{\lambda_j} +\frac{k}{b})| \right) >0.
\ene

Then $\{D_{\lambda_j}T_{kb} \psi \}_{j,k\in\mz}$ is a frame for $\ltr,$ with bounds $A,B.$
\el  Note that Lemma \ref{r-1} contains the case of a standard
wavelet system $\{D_{a^j}T_{kb} \psi \}_{j,k\in\mz}$ as a special case; see
Theorem 15.2.3 in \cite{CB}.

\section{Information packets} \label{251701a}

The purpose of this section is to give a formal introduction to information packets and relate the concept to frames and fusion frames. In particular we will show that information packets generalize as well frames as fusion frames. The definition is as follows.

\bd \label{251701d}  A set  $\{W_{j}\}_{j\in J}$ of closed subspaces of a Hilbert space $\h$ is called an information packet if each $f\in \h$ has an unconditionally convergent expansion
\bee \label{241908a} f= \sum_{j\in J} f_j, \, f_j\in W_j.\ene \ed

\begin{rem} \label{251801a} {\rm If
$\{W_{j}\}_{j\in J}$ is an information packet
and the closed subspaces $U_j, j\in J$ satisfy that
$W_j \subseteq U_j$ for all $j\in J,$ then also $\{U_{j}\}_{j\in J}$ is  an information packet.} \ep \end{rem}

\bex \label{251801b} {\rm Clearly, each frame (or unconditional Schauder basis) $\{f_k\}_{k\in I}$ yields an information packet, with
     one-dimensional spaces $W_j= \Span(f_j), \, \, j\in I.$ More generally, given any
    frame $\{f_k\}_{k\in I}.$  consider
    any collection of sets $\{\sigma_j\}_{j\in J} \subset I$ such that
    \bes \bigcup_{j\in J} \sigma_j=I,\ens
    and let
    \bes W_j:= \span\{f_k\}_{k\in \sigma_j}, \, j\in J.\ens
    Then $\{W_{j}\}_{j\in J}$ is an information packet.}  \ep \enx

In order to relate information packets and fusion frames we need the following lemma.

\bl \label{251801g} Assume that $\{W_j\}_{j\in J}$ is an information packet and that $T: \h \to \h$ is a bounded bijective operator. Then also $\{T(W_j)\}_{j\in J}$ is an information packet. \el

\bp Clearly the subspaces $T(W_j)$ are closed. Furthermore,
given $f\in \h$ we can write $T^{-1}f= \sum_{j\in J} g_j$ for
some $g_j\in W_j,$ with unconditional convergence of the series. Thus, with $f_j:= Tg_j \in T(W_j),$
\bes f= \sum_{j\in J} Tg_j = \sum_{j\in J} f_j,\ens again with unconditional convergence, as desired. \ep

\bl \label{251801f} Let $\{W_j, \omega_j\}_{j\in J}$ be a fusion frame with frame operator $S.$ Then
$\{W_j\}_{j\in J}$ is an information packet. \el

\bp By \eqref{241901e} we clearly have that $\{S^{-1}(W_j)\}_{j\in J}$ is an information packet. Now
the result follows from Lemma \ref{251801g} by letting
$T:=S.$ \ep

The next example shows that information packets is indeed a more general concept than fusion frames, by presenting an information packet $\{ W_j\}_{j\in J}$ for which
$\{ W_j, \omega_j\}_{j\in J}$ is not a fusion frame for any choice of the weights $\omega_j.$

\bex {\rm Let $\etk$ be an orthonormal basis for $\h,$ and let
\bes W_j:= \Span \{e_1, e_j\}, \, j\in \mn.\ens
Applying Example \ref{251801b} and Remark \ref{251801a} to the frame (ONB) $\etk,$  we directly see that
$\{W_{j}\}_{j=1}^\infty $ is an information packet.
However, regardless of the choice of coefficients $\omega_j>0,$
$\{ W_j, \omega_j\}_{j=1}^\infty$ is not a fusion frame. Indeed,
fixing any $k\ge 2,$ we have that
\bes \sum_{j=1}^\infty \omega_j^2 || P_j e_k || ^2
= \omega_k^2.\ens  Thus, if $\{ W_j, \omega_j\}_{j=1}^\infty$ were a fusion frame with bounds $A,B>0,$ then $A \le \omega_k^2 \le B$ for all $k\ge 2.$ On the other hand,

\bes  \sum_{j=1}^\infty \omega_j^2 || P_j e_1 || ^2
= \sum_{j=1}^\infty \omega_j^2,\ens so the fusion property would also imply that $A \le \sum_{j=1}^\infty \omega_j^2 \le B;$ this is clearly not possible, so we conclude that
$\{ W_j, \omega_j\}_{j=1}^\infty$ is not a fusion frame.  } \ep \enx

\section{Weaving of information packets} \label{251901s}

Weaving of information packets is formally defined in the same way as for frames:

\bd \label{241908c} Fix $N\in \mn.$ Then $N$  information packets
$\{W_j^{(\ell)}\}_{j\in J}, \, \ell=1, \dots, N,$ are said to be woven if for each partition $\{\sigma_k\}_{k=1}^N$ of $J,$
\bee \label{241908b} \{W_j^{(1)}\}_{j\in \sigma_1}\cup
\{W_j^{(2)}\}_{j\in \sigma_2} \cup \cdots \{W_j^{(N)}\}_{j\in \sigma_N}
\ene
is an information packet.
\ed

\begin{rem} \label{251901q} That  $N$ information packets
    $\{W_j^{(\ell)}\}_{j\in J}, \, \ell=1, \dots,N,$ are woven means precisely that if we for each $j\in J$ choose
    any $\ell_j\in \{1, \dots, N\},$ then $\{W_j^{(\ell_j)}\}_{j\in J}$ is an information packet.
    \end{rem}

We will now derive a concrete construction of woven information packets, based on certain
wavelet frame constructions in $\ltr.$ We will be dealing with wavelet systems generated by a function $\psi\in \ltr$ such that its Fourier transform $\hp$ is compactly supported;  in the entire discussion we let $I \subset \mr$ denote a finite interval such that
\bes \mbox{supp}\, \hp \subseteq I.\ens

Our main purpose is to prove the following result.

\bt  \label{251701f} Let $\psi \in \ltr.$ Assume that $\hp$ is bounded, compactly supported and that there exist constants $C,D>0$ and
$\alpha, \beta >0$ such that for a neighborhood $U$ of
$0\in \mr,$
\bee \label{232412ba} C\, | \ga|^\beta \le | \hp(\ga)| \le D\, | \ga|^\alpha, \,  \ga \in U. \ene

Fix any $a>1,$ and any $b\in ]0, |I|^{-1}],$  and let

\bee \label{240908ag} V_j= \span \{D_{a^j} T_ {kb}\psi\}_{k\in \mz}, \, j\in \mz.\ene

Fix any $N\in \mn.$ Then   the $N$ collections of subspaces
\bes \{V_{Nj}\}_{j\in \mz},  \{V_{Nj+1}\}_{j\in \mz}, \dots, \{V_{Nj+N-1}\}_{j\in \mz}\ens are woven information packets in $\ltr.$
\et

We will derive Theorem \ref{251701f} as a consequence of the following result.

\bpr \label{240908b}  Let $\psi \in \ltr.$ Assume that $\hp$ is bounded, compactly supported and that there exist constants $C,D>0$ and
$\alpha, \beta >0$ such that for a neighborhood $U$ of
$0\in \mr,$
\bee \label{232412b} C\, | \ga|^\beta \le | \hp(\ga)| \le D\, | \ga|^\alpha, \,  \ga \in U. \ene
Now, fix $a>1$, $ b\in ]0, |I|^{-1}]$ and $N\in\mn$.
Then for any choice of  $\ell_j\in\{0,1,\dots,N-1\}$, $j\in\mz$,
the set $\{ D_{a^{\ell_j}} D_{a^{Nj}}T_{kb}\psi\}_{j,k\in \mz}$ is a frame for $\ltr.$
\epr

\bp We will check the conditions
\eqref{r-1a} and \eqref{r-1b} with
$\lambda_j:=a^{\ell_j+Nj}.$
By assumption, if $0<b \le |I|^{-1},$ then
\bes \sum_{k\neq 0} \sujz |\hp(\frac{\ga}{a^{\ell_j+Nj}})
\hp(\frac{\ga}{a^{\ell_j+Nj}} +\frac{k}{b})|=0\ens for a.e. $\ga\in \mr$ and all $a>1.$
Now choose $J\in \mz$ such that
$([-a^{NJ},-a^{N(J-1)}] \cup [a^{N(J-1)}, a^{NJ}]) \subset U.$
Now, let  $\ga\in \mr \setminus \{0\}$. Then
$|\ga| \in [a^{N(j_0-1)}, a^{Nj_0}]$ for some $j_0\in\mz$.
Choose $K>0$ such that
\bes
&& a^{-N(-J+j_0-K)}\left([- a^{Nj_0},-a^{N(j_0-1)}] \cup [a^{N(j_0-1)}, a^{Nj_0}] \right)\cap I = \emptyset, \ens i.e., 
\bes \left([-a^{N(J+K)},- a^{N(J-1+K)}] \cup [a^{N(J-1+K)},
a^{N(J+K)}] \right) \cap I= \emptyset.
\ens
Then for $j\leq -J+j_0-K-1$, we have
$$a^{-(\ell_j+Nj)}\geq a^{-(N +N(-J+j_0-K-1))}=a^{-N(-J+j_0-K)};$$
thus $\hp(a^{-(\ell_j+Nj)}\ga)=0.$
Since
\begin{eqnarray*}
&&a^{-N(-J+j_0)}\left( [-a^{Nj_0},-a^{N(j_0-1)}] \cup [a^{N(j_0-1)}, a^{Nj_0}]\right) \\
&&=\left([-a^{NJ},-a^{N(J-1)}]\cup [a^{N(J-1)}, a^{NJ}]\right)
\subset U,
\end{eqnarray*}  we have \bes \sujz |
\hp(\frac{\ga}{a^{\ell_j+Nj}})|^2
& = &  \sum_{j=-J+j_0-K}^\infty  | \hp(\frac{\ga}{a^{\ell_j+Nj}})|^2 \\
& = &    \sum_{j=-J+j_0}^\infty  | \hp(\frac{\ga}{a^{\ell_j+Nj}})|^2
+\sum_{j=-J+j_0-K}^{-J+j_0-1}  | \hp(\frac{\ga}{a^{\ell_j+Nj}})|^2 \\
& \le &   \sum_{j=-J+j_0}^\infty
D^2 a^{2\alpha(-\ell_j-Nj)} |\ga|^{2\alpha} +K\, || \hp||_\infty^2\\
& \le &   \sum_{j=-J+j_0}^\infty
D^2 a^{-2\alpha Nj} |\ga|^{2\alpha} +K\, || \hp||_\infty^2\\
& = & D^2 \frac{a^{-2\alpha N(-J+j_0)}}{1- a^{-2\alpha N}}a^{2\alpha Nj_0}+K\, || \hp||_\infty^2\\
& = & D^2 \frac{a^{2\alpha NJ}}{1- a^{-2\alpha N}}+K\, || \hp||_\infty^2.
\ens
Now, still for
$|\ga| \in [a^{N(j_0-1)}, a^{Nj_0}]$ for some  $j_0\in\mz$,
\begin{eqnarray*}
    \sujz | \hp(\frac{\ga}{a^{\ell_j+Nj}})|^2
    &\ge&  | \hp(\frac{\ga}{a^{\ell_j+N(-J+j_0)}})|^2 \\
    &\ge&
    \frac{C^2 |\ga|^{2\beta}}{a^{2\beta(\ell_j+N(-J+j_0))}}
\\  &\ge&
    \frac{C^2 a^{2\beta N(j_0-1)}}{a^{2\beta(N-1+N(-J+j_0))}}
\\  &\ge&
    \frac{C^2  a^{-2\beta N }}{ a^{2\beta (N-1- NJ) }}.
\end{eqnarray*}
Thus  the  frame conditions \eqref{r-1a} and \eqref{r-1b} in Lemma \ref{r-1}
is satisfied, i.e., \\ $\{ D_{a^{\ell_j}} D_{a^{Nj}}T_{kb}\psi\}_{j,k\in \mz}$ is a frame for $\ltr.$
\ep

\vn{\bf Proof of  Theorem \ref{251701f}:} In the entire proof
we fix the assumptions and notation from Theorem \ref{251701f}.
By Proposition \ref{240908b} $ \{D_{a^{Nj}}T_{kb}\psi\}_{j,k\in \mz}$ is a frame; using the associated frame decomposition
of $f\in \ltr$ (here $S$ denotes the frame operator),
\bes f= \sum_{j\in \mz} \sum_{k\in \mz} \la f, \si
D_{a^{Nj}}T_{kb}\psi \ra D_{a^{Nj}}T_{kb}\psi, \ens
it directly follows that the collection of spaces $ \span\{V_{Nj}\}_{j\in \mz}$ is an information packet. Similarly, fixing any $\ell\in \{0,1, \dots, N-1\},$
$ \{ D_{a^{\ell}} D_{a^{Nj}}T_{kb}\psi\}_{j,k\in \mz}$ is a frame, i.e.,
$ \{  D_{a^{Nj+\ell}}T_{kb}\psi\}_{j,k\in \mz}$
is a frame; this implies  that each of the collections
\bee \label{251901p} \{V_{Nj}\}_{j\in \mz},  \{V_{Nj+1}\}_{j\in \mz}, \dots, \{V_{Nj+N-1}\}_{j\in \mz}\ene are information packets in $\ltr.$

However, Proposition \ref{240908b}  yields an even stronger conclusion, namely, that for
any choice of  $\ell_j\in\{0,1,\dots,N-1\}$, $j\in\mz$,
the set $\{ D_{a^{\ell_j}} D_{a^{Nj}}T_{kb}\psi\}_{j,k\in \mz}$ is a frame. This means that for each $\ell_j\in \{0,1, \dots, N-1\},\, j\in \mz,$
the collection of spaces $\{V_{Nj+\ell_j}\}_{j\in \mz}$ is an information packet, i.e., that the information packets in \eqref{251901p} are woven. \ep

The weaving procedure specified in Theorem \ref{251701f} can  be transferred to Gabor systems in a weaker form. The proof is very similar to the proof of  Theorem \ref{251701f}, so we will only sketch it.

\bex  \label{252001b} {\rm
Like wavelet systems, Gabor systems in $\ltr$ are generated by two parameters; indeed, fixing $a,b>0$ and $g\in \ltr,$ a Gabor system has the form
        $\mts$ where
    $E_{mb}T_{na}g(x):= e^{2\pi i mbx}g(x-na).$
    Recall that
    a necessary condition for such a family to be a frame
    is that $ab \le 1.$

        Now fix $a>0$ and $N\in \mn, $ and choose a bounded compactly supported function $g\in \ltr$ such that
        $|g(x)| \ge \epsilon>0$ for all $x\in [0, aN].$
        Letting $I$ denote a finite interval containing the support of $g,$ and choose $b>0$ such that $ b\le |I|^{-1}.$ These assumptions immediately imply that
        \bes \sum_{k\in \mz \setminus \{0\}} \left|\sum_{n\in \mz} g(x-nNa) \overline{g(x-nNa - k/b)} \right|=0,\ens and that, for some constant $C>0,$
            \bes  \epsilon^2 \le \sum_{n\in \mz} |g(x-nNa)|^2 \le C, \, x\in \mr.\ens
    By Theorem 11.4.2 in \cite{CB} this implies that
    $\{E_{mb}T_{naN} g\}_{m,n\in \mz}$
    is a frame for $\ltr.$  As a consequence, fixing any $\ell \in \{0, 1, \dots, N-1\},$  the system of functions
        $\{T_{\ell a}E_{mb}T_{naN} g\}_{m,n\in \mz}$
    is a frame; using the commutator relations between translation and modulation operators, this implies that
    $\{E_{mb}T_{(nN + \ell) a} g\}_{m,n\in \mz}$
    is a frame. This conclusion can be phrased in terms of information packets. Indeed, for $n\in \mz,$ let
    \bee \label{252001a}  V_n:= \span\{ E_{mb}T_{n  a}g\}_{m\in \mz};\ene
    then each of the collections of spaces
    \bee \label{251901pa} \{V_{Nn}\}_{n\in \mz},  \{V_{Nn+1}\}_{n\in \mz}, \dots, \{V_{Nn+N-1}\}_{n\in \mz}\ene are information packets in $\ltr.$ Compare this conclusion with \eqref{251901p}!

    As for the wavelet case, the weaving property for the information packets in
    \eqref{252001a} now follows from the standard technical condition for irregular Gabor systems \cite{CC5}. Assume that we for each $n\in \mz$ choose $\ell_n \in \{0, 1, \dots, N-1\}.$ Letting now $\ga_n:= nNa+ \ell_n a,$ similar considerations as above yield that
        \bes \sum_{k\in \mz \setminus \{0\}} \left|\sum_{n\in \mz} g(x-\ga_n)   \overline{g(x-\ga_n - k/b)} \right|=0,\ens and that, for some constant $C>0,$
        \bes  \epsilon^2 \le \sum_{n\in \mz} |g(x-\ga_n)|^2 \le C, \, x\in \mr.\ens

By   Corollary 2.2 in \cite{CC5}, this implies that
    $\{E_{mb}T_{\ga_n a} g\}_{m,n\in \mz} =\{E_{mb}T_{(nNa+ \ell_n a) a} g\}_{m,n\in \mz} $
    is a frame. This means that for each $\ell_n\in \{0,1, \dots, N-1\},\, n\in \mz,$
    the collection of spaces $\{V_{Nn+\ell_n}\}_{n\in \mz}$ is an information packet, i.e., that the information packets in \eqref{252001a} are woven.

     The difference to the wavelet case is that
      for a given Gabor frame $\mts,$
      this procedure can only work for  $N\in \mn$
     chosen such that $abN\le 1,$ not for all $N\in \mn.$
} \ep \enx

\begin{rem} {\rm Although the procedures of obtaining woven information packets in Theorem \ref{251701f} and Example \ref{252001b} are very similar,
there are a number of noticeable differences between the wavelet case and the Gabor case:

\bei \item[(i)] In the wavelet case, Theorem \ref{251701f}
yields a construction of a single wavelet system
$\{ D_{a^{j}}T_{kb}\psi\}_{j,k\in \mz},$  i.e., a fixed choice of some $a>1, b>0,$ which can be used to construct $N$ woven information packets for any choice of $N\in \mn.$ In contrast, in Example \ref{252001b}, if we fix
the parameter $a$ in a Gabor system and want
to split the Gabor system into $N$ information packets following the procedure above,  we
will have to choose the parameter $b>0$ such that  $abN \le 1.$
\item[(ii)] In the Gabor case in Example \ref{252001b}, the weaving is taking place
via a partition of the set of  translations $na, \, n\in \mz.$  A completely similar result can be obtained, with weaving taking place via a
partition of the modulations $mb,\, m\in \mz,$ simply by applying Fourier transform techniques (applying the Fourier transform turns translation to modulation and modulation to translation).
In the wavelet case in Theorem \ref{251701f}, the weaving is taking place via a partition of
the scalings $a^j, \, j\in \mz.$ For wavelet systems it is not clear whether it is possible to obtain woven information packets via a
partition of the set of translates $kb, \, k\in \mz.$
\eni \ep}
\end{rem}

\noindent{\bf Acknowledgment:} The authors would like to thank the reviewers for many good remarks, which improved the presentation of the results. They also thank Peter Massopust for discussions and for suggesting the wording {\it information packet} instead of our original {\it information package.}


{\bf \vspace{.1in} \noindent Ole Christensen\\
    Department of Applied Mathematics and Computer Science\\
    Technical University of Denmark,
    Building 303,
    2800 Lyngby, Denmark\\
    Email: ochr@dtu.dk

    \vspace{.1in}\noindent Hong Oh Kim \\
    Department of Mathematical Sciences, KAIST\\
    291 Daehak-ro, Yuseong-gu, Daejeon 34141,
    Republic of Korea\\
    Email: kimhong@kaist.edu

    \vspace{.1in} \noindent Rae Young Kim \\
    Department of Mathematics, Yeungnam University\\
    280 Daehak-Ro, Gyeongsan, Gyeongbuk 38541,
    Republic of Korea\\
    Email:  rykim@ynu.ac.kr}

\end{document}